\documentclass[12,letterpaper]{article}

\usepackage{blindtext} 
\usepackage[sc]{mathpazo} 
\usepackage[OT1]{fontenc} 
\usepackage{microtype} 
\usepackage[english]{babel} 
\usepackage[margin=1.1in,columnsep=18pt]{geometry}  
\usepackage[hang, small,labelfont=bf,up,textfont=it,up]{caption} 
\usepackage{booktabs} 

\usepackage{lettrine} 
\usepackage{graphicx}
\usepackage{enumitem} 
\setlist[itemize]{noitemsep} 

\usepackage{amsmath}  
\usepackage{listings} 
\usepackage{hyperref} 

\usepackage{textcomp}

\usepackage{lipsum}
\usepackage{titlesec}

\usepackage{abstract} 

\usepackage{titlesec} 
\titleformat{\section}{\bfseries\normalsize}{\thesection.}{1em}{} 
\titleformat{\subsection}[block]{\itshape\normalsize}{\thesubsection.}{1em}{} 

\usepackage{fancyhdr} 
\usepackage{titling} 

\usepackage{hyperref} 

\pretitle{\begin{center}\Large\bfseries} 
\posttitle{\end{center}} 
\title{A Solution Algorithm for Minimum Convex Cost Network Flow Problems} 
\author{%
Dewan Ferdous Wahid\thanks{\textit{Corresponding author}}\thanks{Department of Natural Sciences, Stamford University Bangladesh, Email: dfwahid@stamforduniversity.edu.bd}
\qquad 
Farjana Habiba\thanks{Department of Natural Sciences, Stamford University Bangladesh, Email: farjanahabiba4@gmail.com} 
\qquad 
Ganesh Chandra Ray\thanks{Department of Mathematics, University of Chittagong, Email: ganesh.ray@cu.ac.bd
\newline
Copyright \textcopyright $\:$ 2012, Authors}
}
\date{\vspace{-5ex}} 

\begin{document}
\maketitle
\begin{abstract}
This paper presents a heuristic algorithm to solve Minimum Convex-Cost Network Flow Problems (MC-CNFP) heuristically. This solution algorithm is constructed on the concepts of Network Simplex Method (NSM) for minimum cost network flow problem, Convex Simplex Method (CSM) of Zangwill, the decomposition of convex simplex method and non-linear transformation problem. 
\end{abstract}

\section{Introduction}
The minimum convex-cost network flow problem (MC-CNFP) is a class of minimum cost network flow problems with convex cost function. This problem structure may appears in different real-world problems such that cost of power interruption in electrical supply networks, delay cost of communication networks and over-crowding costs in city transportation networks etc.

Consider $G(N, A)$ is a directed network, where $N=\{1,..,m\}$ and $A=\{(i, j),..,(s, t)\}\subseteq N×N$ are node and arcs sets respectively. Let $x_{ij}$ be the flow through the arc$(i j)$, and the vector $x = \{x_{ij}: \: (i,j) \in  A)\}$. Then MC-CNFP can be formulated as
\begin{flalign}
&\min \: \sum\sum_{(i,j)\in A} c_{ij}(x_{ij}),\\
&\text{subject to:} \:\sum_{j:\: (i,j)\in A}x_{ij} - \sum_{k:\: (k,i)\in A}x_{ki} = b_i; \qquad \forall i \in N,\\
&\qquad\qquad\qquad x_{ij}\geq 0; \qquad \forall (i,j)\in A,
\end{flalign}
where $b_i$ is the net flow generated at node $i$ and $c_{ij}: \: \mathbb{R} \rightarrow \mathbb{R}$ are given convex cost functions with continuous first derivative for arcs $(i, j)$. The above formulation also written as
\begin{flalign}
&\min \: C(x), \label{eq4}\\
&\text{subject to:} \: Ax = b; \qquad x\geq 0,
\end{flalign}
where $C(x)$ is convex and the constraints are linear equations. The matrix $A$ is the node-arc incidence matrix with rank $(m−1)$ \cite{bazaraa2005linear}.
\\
\\
This paper represents an optimality condition to minimize the objective function in Eq.(\ref{eq4})  with subject to linear constraints.

\section{Condition for Optimality}
We introduce an artificial arc to root node (any other node would do), that lead to the extended constraint matrix $A_e = (A, e_m)$ of rank $m$, where $e_m$ is a unit vector \cite{bazaraa2005linear}. 

Then Eq.(\ref{eq4}) can rewrite as
\begin{flalign}
&\min \:\: C(x_e), \label{eq6}\\
&\text{subject to:}\: \: A_e x_e = b; \label{eq7}\\
& x_e \geq 0, \label{eq8}
\end{flalign}
where $x_e$ is $n\times 1$ and $A_e$ is $m\times 1$ matrix, here $n$ is the number of are including with artificial arc. Now, the Lagrangian for Eq.(\ref{eq6}) can be formulated as
\begin{flalign}
z(x_e, \mu, \lambda) = C(x_e) + \mu^T (b - A_e x_e) - \lambda x_e,
\end{flalign}
where $\lambda$ and $\mu$ are Lagrange multipliers. The optimum value $\bar{x}$ of Eq.(\ref{eq6}) should satisfy the following Karush-Kuhn-Tucker (KKT) conditions \cite{zangwill1967convex}.
\begin{flalign}
&\nabla z = \nabla C(\bar{x}) - \mu^T A_e - \lambda = 0,\\
&\lambda \bar{x} = 0,\\
&\bar{x} \geq 0,\\
&\lambda \geq 0.\\
\end{flalign}
For each arc flow $x_{ij}$ associated with the arc$(i,j)$, we get
\begin{flalign}
&\frac{\delta z}{\delta x_{ij}} = \frac{\delta C(x)}{\delta x_{ij}} - \mu^T a_{ij} - \lambda_{ij} = 0, \label{eq14}\\
&\lambda_{ij}x_{ij} = 0,\\
&x_{ij} \geq 0,\\
&\lambda_{ij}\geq 0 \label{eq17}.
\end{flalign}
where $\mu^T \in \mathbb{R}_n$ and $a_{ij}$ is column vector associated to $x_{ij}$ (has value $1$ in the $i$-th row and $-1$ in the $j$-th row in $A_e$). Therefore, from Eq.s(\ref{eq14})-(\ref{eq17}) we get
\begin{flalign}
&\frac{\delta z}{\delta x_{ij}} = \frac{\delta C(x)}{\delta x_{ij}} - (\mu_i -\mu_j) - \lambda_{ij} = 0, \label{eq18}\\
&\lambda_{ij}x_{ij} = 0,\\
&x_{ij} \geq 0,\\
&\lambda_{ij}\geq 0 \label{eq21}.
\end{flalign}
By using the Eq.s(\ref{eq18})-(\ref{eq21}), we get
\begin{flalign}
&\frac{\delta z}{\delta x_{ij}} = \frac{\delta C(x)}{\delta x_{ij}} - (\mu_i -\mu_j) \geq 0, \label{eq22}\\
&x_{ij}\frac{\delta z}{\delta x_{ij}} = x_{ij}\big[\frac{\delta C(x)}{\delta x_{ij}} - (\mu_i -\mu_j)\big] = 0,\\
&x_{ij} \geq 0, \label{eq24}
\end{flalign}
Therefore, $\bar{x}$ will minimize the MC-CNFP, in Eq.(\ref{eq6}), if it satisfies the optimality conditions given in Eq.s(\ref{eq22})-(\ref{eq24}).

\section{Solution Algorithm}
In the first step of the algorithm, we find an initial basic feasible solution. Then, we use iterative procedure for moving towards optimal solution.

\subsection{Initial Feasible Solution}
Since, the constraints, in Eq.s(\ref{eq7})-(\ref{eq8}), are linear, we follow the procedure of the inspection of a spanning tree (basis sub-grap) similar to the network simplex method with linear constraints discussed in Bazarra et al. Let $\bar{x}^0 = (\bar{x}^0_B, \bar{x}^0_N)$ is a initial feasible solution, where $\bar{x}^0_B$ and  $\bar{x}^0_N$ are the basic and nonbasic solutions respectively. In next, we have to improve this initial feasible solution to approximate an optimal solution.

\subsection{Testing Optimality of a Feasible Solution}
Any feasible solution of Eq.(\ref{eq6}) would be the optimal solution if it satisfies the conditions in Eq.s(\ref{eq22})-(\ref{eq24}). Let $\bar{x}^k = (\bar{x}^k_B, \bar{x}^k_N)$ be a feasible solution in any $k$-th iteration and $I^k_B = \{ij : \:\: x_{ij}^k \in \bar{x}^k_B\}$ , $I^k_N = \{ij : \:\: x_{ij}^k \in \bar{x}^k_N\}$. We have $x_{ij}^k \ge 0$, then the complementary slackness condition implies that 
\begin{flalign}
&\frac{\delta z}{\delta x_{ij}} = \frac{\delta C(x)}{\delta x_{ij}} - (\mu_i -\mu_j) = 0; \: \forall ij \in I^k_B.
\end{flalign}
Let,
\begin{flalign}
&\frac{\delta z}{\delta x_{rl}^k} = \min \big\{ \frac{\delta z}{\delta x_{ij}^k}; \: ij \in I^k_N\big\},\\
&x_{st}^k \: \frac{\delta z}{\delta x_{st}^k} = \max \big\{ \frac{\delta z}{\delta x_{ij}^k}; \: ij \in I^k_N\big\},
\end{flalign}

\label{them-1}\noindent\textbf{Theorem 3.1}:	\textit{If $\big|\frac{\delta z}{\delta x_{rl}^k}\big| = x_{st}^k \: \frac{\delta z}{\delta x_{st}}$, for a feasible solution $\bar{x}$, then $\bar{x}$ is optimal.} \\\\

\noindent\textit{Proof}: Let $\big|\frac{\delta z}{\delta x_{rl}^k}\big| = x_{st}^k \: \frac{\delta z}{\delta x_{st}}$. Then we have, 
\begin{flalign}
&\frac{\delta z}{\delta x_{ij}} = \frac{\delta C(x)}{\delta x_{ij}} - (\mu_i -\mu_j) = 0, \qquad \text{if } x_{ij} \ge 0, \label{eq9}\\
&\frac{\delta z}{\delta x_{ij}} = \frac{\delta C(x)}{\delta x_{ij}} - (\mu_i -\mu_j) \ge 0, \qquad \text{if } x_{ij}=0. \label{eq10}
\end{flalign}
Here Eq.s(\ref{eq9})-(\ref{eq10}) and the feasibility of $\bar{x}$ are simply the conditions in Eq.s(\ref{eq22})-(\ref{eq24}), which also provides the condition for optimality for the problem given in Eq.s(\ref{eq6})-(\ref{eq7}) \cite{hsia1973decomposition}.
  
\subsection{Iterative Procedure for Moving Towards Optimal Solution}
Any feasible solution which fails to satisfy the optimal condition, in \hyperref[them-1]{\textit{Theorem 3.1}}, has to improve toward the optimal solution by changing nonbasic variables to basic. Since the objective function in Eq.(\ref{eq6}) is convex, so here we use an iterative procedure described in \cite{hsia1973decomposition} for convex simplex method. To improve a feasible solution following cases need to be considered;
\\\\
\label{c1}\noindent\textit{Case 1}: If  $\big|\frac{\delta z}{\delta x_{rl}^k}\big| \ge x_{st}^k \: \frac{\delta z}{\delta x_{st}}$; increase $x_{rl}^k$ by $\Delta^k$, where $\Delta^k$ is compute as following procedure.

Let $I_{B_{rl}}^k = \{lu,..., ij, ..., wr\}$ be the set of indices of the basic flows of the loop contacting the arc$(r,l)$ according to the loop direction. Then 
\begin{flalign}
\Delta^k = \min \big\{|x_{ij}| : \: ij \in I^k_{B_{rl}} \text{ and } x_{ij} \in \bar{x}^k_B \big\}.
\end{flalign}
				
\begin{figure}[h!] \label{f1}
	\centering
	\includegraphics[width=80mm]{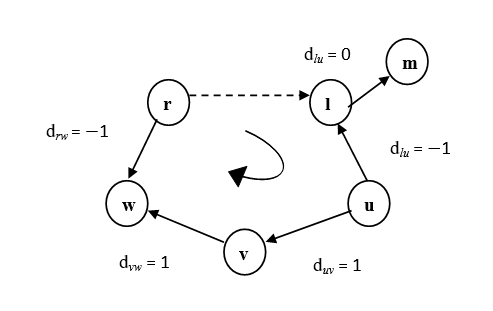}
	\caption{Direction of the basic loop containing the arc$(r,l)$.}
\end{figure}
Next adjust the flow of the network according to loop direction (\hyperref[f1]{Fig. 1}) as follows
\begin{flalign}
&y_{ij}^k = x_{ij}^k; \: ij \in I_N^k -\{rl\},\\
&y_{rl}^k = x_{rl}^k + \Delta^k,\\
&y_{ij}^k = x_{ij}^k + d_{ij}\Delta^k; \: \forall \: ij \in I_B^k, 
\end{flalign}
where \begin{flalign}
d_{ij} = \begin{cases} 
1;  & ij \in I_{B_{rl}}^k, \\
-1; & ji \in I_{B_{rl}}^k, \\
0;  & ij, ji \in I_{B_{rl}}^k, 
\end{cases}\label{eq35}
\end{flalign}
By doing so, one the basic flow say $x_{B_{ij}}^k$ may be driven to zero. Let $\bar{y}^k$ be the value of $\bar{x}^k$ after making the necessary adjustment. Since the function is convex, so a better point could be found before reaching $\bar{y}^k$ \cite{bazaraa2006nonlinear}.

To check this, find $\bar{x}^{k+1}$ by using the line search
\begin{flalign}\label{ls}
C(\bar{x}^{k+1}) = \min \{C(\bar{x}^k); \: \bar{x} = \lambda \bar{x}^k + (1-\lambda)\bar{y}^k, \: 0 < \lambda < 1 \}
\end{flalign}
If $\bar{x}^{k+1}\neq \bar{y}^k$ do not change the former basis and go to the next iteration. If $\bar{x}^{k+1}= \bar{y}^k$ and if a basic flow becomes zero during the adjustment made, change the former basic and go to the next iteration.\\
\\
\label{c2}\noindent\textit{Case 2}: If  $\big|\frac{\delta z}{\delta x_{rl}^k}\big| < x_{st}^k \: \frac{\delta z}{\delta x_{st}}$; decrease $x_{st}^k$ by $\Delta^k$, where $\Delta^k$ is determined as in \hyperref[c2]{\textit{Case 1}}.

Next adjust the flow of the network as follows
\begin{flalign}
&y_{ij}^k = x_{ij}^k; \: ij \in I_N^k -\{st\},\\
&y_{rl}^k = x_{rl}^k - \Delta^k,\\
&y_{ij}^k = x_{ij}^k - d_{ij}\Delta^k; \: \forall \: ij \in I_B^k, 
\end{flalign}
where $d_{ij}$ can calculate as in Eq.\ref{eq35}.

Then we obtain $\bar{y}^k$. As we decrease $x_{st}^k$, then either $x_{st}^k$ itself or any basic flow say will be driven to zero. Now calculate $\bar{x}^{k+1}$ from the line search in Eq.(\ref{ls}). If $\bar{x}^{k+1} \neq \bar{y}^k$, do not change the former basis and go to next iteration and if $\bar{x}^{k+1} = \bar{y}^k$ change the basis. 

\subsection{Optimality Condition During Line Search}
In line search problem, we search toward the optimal solution by solving 
\begin{flalign}\label{ls-1}
C(\bar{x}^{k+1}) = \min \{C(\bar{x}^k); \: \bar{x} = \lambda \bar{x}^k + (1-\lambda)\bar{y}^k, \: 0 < \lambda < 1 \},
\end{flalign}
where $\bar{x}^{k+1} = \lambda \bar{x}^k + (1-\lambda)\bar{y}^k$. However, from practical experience, in the case of some problem, we see $\lambda=1$. Therefore in this case, $\bar{x}^{k+1} = \bar{x}^k$, i.e. the line search problem indicates that there is no other better point close to the optimal solution than the point $\bar{x}^k$. Again, if we proceed to the next iteration then the feasible solution will not change and the problem will circulate to the point  $\bar{x}^k$ with satisfying the optimality condition given in \hyperref[them-1]{\textit{Theorem 3.1}}. But this feasible solution makes the cost function least compared to feasible solutions in the previous iterations. Therefore, to avoid the circular situation, we return the point $\bar{x}^k$ as an heuristic solution if $\lambda=1$.
\section{Conclusion}
In this paper, we propose a heuristic solution procedure for minimum convex-cost network flow problem. From implementation experience, we set another optimality condition without giving any strong mathematical logic to avoid the circular situation when $\lambda=1$. Therefore, the implementation for large scale data and demonstrating any logical  explanation for the special condition are still open.  

\small
\bibliographystyle{IEEEtran}
\bibliography{mccnfp-cu}
\end{document}